\documentclass[12pt]{amsart}
\usepackage[margin=1in]{geometry}
\usepackage{amssymb, amsthm, amsmath}
\usepackage[mathscr]{euscript}
\usepackage{graphicx, xcolor}
\usepackage[colorlinks,bookmarks]{hyperref}
\hypersetup{%
  linkcolor = blue,
  citecolor = purple,
  urlcolor = purple,
}
\usepackage{listings}
\theoremstyle{definition}

\theoremstyle{plain}
\newtheorem{thm}{Theorem}

\theoremstyle{definition}

\newtheorem{eg}{Example}

\theoremstyle{remark}

\newtheorem*{note}{Note}

\newcommand\gl{\mathfrak{gl}}
\newcommand\g{\mathfrak{g}}
\newcommand\h{\mathfrak{h}}
\renewcommand\SS{\mathbf{S}}

\DeclareMathOperator{\ch}{ch}
\DeclareMathOperator{\atyp}{atyp}
\DeclareMathOperator{\typ}{typ}

\begin{document}
\author{Abhik Pal}
\title[Lie Superalgebras in Sage]{%
  Lie Superalgebras and Generalized Kazhdan-Lusztig Polynomials
}
\begin{abstract}
  We present \texttt{liesuperalg} a SageMath package for representation-theoretic calculations involving Lie Superalgebras in Type A.
  Our package introduces functionality to calculate invariants of weights and produce the associated cup diagrams.
  We expose functionality to calculate characters of irreducible representations, work with combinatorics of generalized Kazhdan-Lusztig polynomials, and determine composition factor multiplicities of indecomposable modules.
  Our package implements an algorithm to decompose arbitrary modules in terms of irreducible ones in the Grothendeick group of Lie superalgebra representations.
\end{abstract}
\maketitle
\setcounter{tocdepth}{2}

\section{Introduction}
A fundamental problem in representation theory of Lie algebras is to classify all finite dimensional representations of a given Lie algebra \(\g\) and to determine their characters.
When \(\g\) is a complex semisimple Lie algebra, the theorem of highest weight proves that there is a bijection between equivalence classes of irreducible \(\g\)-modules and \emph{dominant integral weights} where \emph{weights} are generalized eigenvalues with respect to a distinguished subalgebra \(\h\) of \(\g\).
The character formulas in this case are determined by exploiting a relationship between representations of the Lie algebra \(\g\) and that of its underlying \emph{Weyl group}.
This material is fairly standard and an in-depth treatment appears in \cite{FH:RepresentationTheory,Hum:CategoryO,Hum:IntroductionLie,Kna:RepresentationTheory}.

In the 1970s, motivated by mathematical physics, the theory was generalized to Lie superalgebras \cite{cns1975:graded}.
That is, Lie algebras of generalized groups called Lie supergroups typically realized as endomorphism algebras of \(\mathbf Z_2\)-graded vector spaces along with a generalized Lie bracket and Jacobi identity.
The representation theory of such algebras was formalized by Kac in \cite{kac1977:lie} where he also posed the question of calculating the characters of the irreducible modules.
There are several obstructions to generalizing the theory of weights from \emph{ordinary} representation theory to \(\mathbf Z_2\)-graded representation theory: in general the indecomposable and the irreducible modules do not coincide and unlike ordinary representation theory, the representation theory depends on the specific choice of the Borel subgroup \cite{Mus:LieSuperalgebras}.
After a series of several conjectures and partial results (for instance \cite{vhkt1990:character, kw1994:integrable}) the character formula for irreducibles of the general linear Lie superalgebra and composition factor multiplicities were determined by Serganova in the late 1990s \cite{ser1998:characters}.
One of the main ingredients in this result were generalized Kazhdan-Lusztig polynomials \(K_{\lambda,\mu}(q)\) which, when evaluated at \(q = -1\), determine multiplicities appearing in the character formula.

We introduce a SageMath package \texttt{liesuperalg} to help with calculations in the representation theory of the general linear superalgebra.
The package includes functionality to work with weights and calculate invariants that are of interest in the literature.
Additionally, we introduce functionality to produce cup diagrams from weights, these combinatorial gadgets give a convenient way of working with weights of Lie superalgebras and have been used extensively in representation theory \cite{bs2010:highest-1,bs2010:highest-2,bs2010:highest-3,bs2010:highest-4}. The package also computes the generalize Kazhdan-Lusztig polynomials using a version of Brundan's algorithm \cite{bru2003:kazhdanlusztig} as presented in \cite{sz2007:character}.
We expose methods to determine composition factors of indecomposable modules; this is based on the calculation of the Jantzen filtration of such modules in \cite{sz2012:generalised, su2006:composition}.
Our package is also able to decompose arbitrary modules in terms of irreducible characters in the Grothendeick group of \(\gl(m|n)\)-representations.

\subsection{Outline}
We will interleave the description of the package functions with the discussion of the mathematical background.
Most functionality is described through examples computed within Sage sessions. 
In Section~\ref{sec:invariants} we set notation, state the basic facts regarding weights of the general linear superalgebra, introduce the necessary invariants, and describe the diagrams associated to weights.
Representation, character formulas, and generalized Kazhdan-Lusztig polynomials are described in Section~\ref{sec:reps}.
Section~\ref{sec:composition} explains the main decomposition algorithm for Kac modules.
Finally in Section~\ref{sec:decomposition} we describe our approach to decomposing Lie superalgebra modules in terms of irreducible characters. 

\subsection{Supplementary Material}
A up-to-date version of the \texttt{liesuperalg} package is available at the git repository \href{https://github.com/abhikpal/liesuperalg}{\texttt{github.com/abhikpal/liesuperalg}}.
All examples discussed here are viewable in the notebook \texttt{liesuperalg/Examples.ipynb}.

\subsection{Acknowledgments}
We thank Steven Sam for many helpful discussions and feedback.
The author was partially supported by NSF grant DMS 2302149.

\section{Invariants and Diagrams of Weights}
\label{sec:invariants}

We state some of the basic facts related to the general linear Lie superalgebra \(\gl(m|n)\). 
Further details and preliminaries on Lie superalgebras in general appear in Kac's original paper \cite{kac1977:lie} and in textbooks \cite{Mus:LieSuperalgebras,Sch:TheoryLie,CW:DualitiesRepresentations}

\subsection{Basic Definitions}

We fix a \(\mathbf Z_2\)-graded vector space  \(V = V_0 \oplus V_1\) with \(\dim{V_0} = m\) and \(\dim{V_1} = n\). We call the subspaces \(V_0\) and \(V_1\) \emph{even} and \emph{odd} respectively. The \(\mathbf Z_2\)-graded endomorphisms of such a vector space along with a \(\mathbf Z_2\)-graded commutator forms the \emph{general linear Lie superalgebra} \(\g = \gl(m|n)\).

Fix a basis \(E_1, \dots E_{m + n} \) of \(V\) such that \(E_1, \dots, E_m\) form a basis of \(V_0\) and \(E_{m + 1}, \dots, E_{m + n}\) form a basis of \(V_1\).
We work with respect to the Borel subgroup of upper triangular matrices and the Cartan subalgebra \(\h\) of diagonal matrices.
Let \(\epsilon_1, \dots, \epsilon_{m + n}\) be a basis of \(\h^\ast\).
For convenience, set \(\delta_a = \epsilon_{\dot{a}}\) where \(\dot{a} = a + m\).
The space \(\h^\ast\) carries a bilinear form
\begin{equation}
  \label{eq:def:bilinear}
  (\epsilon_i, \epsilon_j) = (-1)^{[i]}\delta_{ij}
\end{equation}
where \([i] = 0\) if \(i \leq m\) and \([i] = 1\) otherwise and \(\delta_{ij}\) is the Kronecker delta.
The sets of \emph{positive even roots} and \emph{odd roots} are respectively given by
\begin{align}
  \Delta_0^+ &=  \left\{ \epsilon_i - \epsilon_j, \delta_k - \delta_l \mid 1 \leq i < j \leq m, 1 \le k < l \leq n \right\}, \\
  \Delta_1^+ &= \left\{ \epsilon_i - \delta_j \mid 1 \leq i \leq m, 1 \leq j \leq n \right\}.
\end{align}

A \emph{weight} of \(\gl(m|n)\) is an element \(\lambda \in \h^\ast\) written as
\begin{equation}
  \label{eq:def:weight}
  \lambda
  = (\lambda^0_1, \dots, \lambda^0_m \mid \lambda^1_1, \dots, \lambda^1_n)
  = \sum_{i = 1}^m \lambda^0_i \epsilon_i - \sum_{j = 1}^n \lambda^1_j \delta_j.
\end{equation}
in the \(\epsilon\delta\)-basis.
We will use the indexing convention 
\begin{equation*}
  \lambda_i =
  \begin{cases}
    \lambda^0_i &i \leq m \\
    \lambda^1_{i - m} & i > m
  \end{cases}.
\end{equation*}

The object \texttt{Weight()} is the main class that implements functionality related to a weight of \(\gl(m|n)\). This is instantiated by providing, as a list, the entries of \(\lambda^0\) and \(\lambda^1\). Several of the parameters related to a weight can be accessed as properties:
\begin{itemize}
\item \texttt{Weight.L} and \texttt{Weight.R} respectively equal \(\lambda^0\) and \(\lambda^1\).
\item \texttt{Weight.coeff\_L} and \texttt{Weight.coeff\_R} compute the coefficients of \(\epsilon_i\) and \(\delta_j\) appearing in \eqref{eq:def:weight}.
\item \texttt{Weight.m} and \texttt{Weight.n} give the dimensions 
  \(m\) and \(n\)
\end{itemize}
The \texttt{Weight} object supports two common indexing conventions:
\begin{itemize}
\item \texttt{Weight[p, i]} produces \(\lambda^p_i\)
\item \texttt{Weight[j]} produces the coefficient of \(\epsilon_j\)  in \eqref{eq:def:weight}
\end{itemize}

A weight is called \emph{integral} when all \(\lambda_i\) are  integers and \emph{dominant} when
\begin{equation}
  \label{eq:def:dominant}
  \lambda^0_1 \geq \dots \geq \lambda^1_m,
  \qquad
  \lambda^1_1 \leq \dots \leq \lambda^1_n.
\end{equation}
We check whether a weight is dominant using the method \texttt{Weight.is\_dominant()}

The Weyl group of \(\gl(m|n)\) is the product of two symmetric groups \(W = \mathfrak{S}_m \times \mathfrak{S}_n\).
This group acts on the set of positive weights by the \emph{dot action} 
\begin{equation}
  \label{eq:def:dot-action}
  \sigma \bullet \mu = \sigma(\mu + \rho) - \rho
\end{equation}
where \(\sigma \in W\), \(\mu\) is a positive weight, and 
\begin{equation}
  \label{eq:def:rho}
  \rho = (m, \dots, 2, 1 \mid 1, 2, \dots, n).
\end{equation}
We will often denote the \emph{\(\rho\)-translate} of \(\lambda\) by
\begin{equation}
  \label{eq:def:rho-translate}
  \lambda^\rho = \lambda + \rho
\end{equation}
The function \texttt{Weight.rho()} calculates the \(\rho\)-translate of any weight and the function \texttt{rho(m, n)} produces the distinguished weight \(\rho\) for \(\gl(m|n)\).
The function \texttt{one(m, n)} produces the distinguished weight \((1, \dots, 1 \mid 1, \dots, 1)\).

Finally, addition, negation, subtraction, and equality testing for weights is defined by overloading the \texttt{+}, \texttt{-}, and \texttt{==} operators.

\begin{eg}
  Let 
  \(\lambda = (7, 6, 5, 5, 3, 3, 2, 2, 0 \mid 1, 2, 3, 4, 4, 5, 7, 7)\) be a weight of \(\gl(9\mid8)\). 

\begin{lstlisting}
sage: from liesuperalg import Weight
sage: w = Weight([7,6,5,5,3,3,2,2,0], [1,2,3,4,4,5,7,7])
sage: w
gl(9|8) weight (7, 6, 5, 5, 3, 3, 2, 2, 0 | 1, 2, 3, 4, 4, 5, 7, 7)
sage: w.L, w.R
([7, 6, 5, 5, 3, 3, 2, 2, 0], [1, 2, 3, 4, 4, 5, 7, 7])
sage: w.m, w.n
(9, 8)
\end{lstlisting}

  Since \(\rho = (9, \dots, 1 \mid 1, \dots, 8)\), the \(\rho\)-translate of \(\lambda\) above is
  \begin{equation*}
    \lambda^\rho = 
    (16, 14, 12, 11,  8,  7,  5,  4,  1 \mid  2,  4,  6,  8,  9, 11, 14, 15).
  \end{equation*}
  Continuing the sage session above, we see
\begin{lstlisting}
sage: from liesuperalg import rho
sage: rho(9, 8)
gl(9|8) weight (9, 8, 7, 6, 5, 4, 3, 2, 1 | 1, 2, 3, 4, 5, 6, 7, 8)
sage: w.rho()
gl(9|8) weight (16, 14, 12, 11, 8, 7, 5, 4, 1 | 2, 4, 6, 8, 9, 11, 14, 15)
sage: w + rho(9, 8) == w.rho()
True
sage: w.rho() - w == rho(9, 8)
True
sage: form liesuperalg import one
sage: rho(5, 4) == one(5, 4)
False
\end{lstlisting}

\end{eg}

\begin{note}
  The class is designed to work with Sage's in-built \texttt{latex()} method and respects the \texttt{\%display typeset} magic command in the notebook interface.
  So calling \texttt{latex()} on a \texttt{Weight} object produces the equivalent \LaTeX\ code:

\begin{lstlisting}
sage: from liesuperalg import rho
sage: latex(rho(5, 4))
\left(5, 4, 3, 2, 1 \mid 1, 2, 3, 4\right)
\end{lstlisting}

\end{note}

\subsection{Representation Theoretic Invariants}
We now list some invariants for weights that play an important role in the representation theory of $\gl(m|n)$.
All of these invariants are implemented as properties of the \texttt{Weight} object using the \texttt{@property} decorator.
This ensures that properties cannot be accidentally overwritten and allows us to only calculate invariants when requested by the user. 

A positive odd root \(\gamma \in \Delta_1^+\) is called an \emph{atypical root} of \(\lambda\) if and only if \((\gamma, \lambda^\rho) = 0\), that is, \(\gamma\) is orthogonal to \(\lambda^\rho\) with respect to the form \eqref{eq:def:bilinear} \cite{sz2007:character, vz1999:characters}.
If \(\gamma = \epsilon_j - \delta_k\), the definition is equivalent to
\begin{equation*}
  (\lambda^\rho, \gamma) =
  (\lambda + \rho, \epsilon_j - \delta_k) = 
  (\lambda_i + m + 1 - j) - (\lambda_{\dot{k}} + k) =
  0.
\end{equation*}
The set of all atypical roots is denoted \(\Gamma(\lambda)\) and the total number of such roots is called the \emph{degree of atypicality} \(r = |\Gamma(\lambda)|\) of \(\lambda\).
We order the atypical roots by the indices of the positive roots. So if \(\Gamma(\lambda) = \left\{ \gamma_1, \dots, \gamma_r \right\}\) we say \(\gamma_1 = \epsilon_{m_1} - \delta_{n_1} < \dots < \gamma_r = \epsilon_{m_r} - \delta_{n_r}\) if \(m_r < \dots < m_1\) and \(n_1 < \dots < n_r\).
We call \(\gamma_s\) the \(s\)-th atypical root of \(\lambda\).
Note that the  indices \(\{ (m_s, n_s) \mid 1 \leq s \leq r \}\) is the maximal set of pairs \((j, k)\) satisfying \(\lambda^\rho_j = \lambda^\rho_{\dot{k}}\).
The \emph{atypical tuple} of \(\lambda\) is defined as
\begin{equation*}
  \atyp(\lambda) = (\lambda^\rho_{\dot{n_1}}, \dots, \lambda^\rho_{\dot{n_r}})
\end{equation*}
and the \emph{typical tuple} is the tuple in \(\mathbf Z^{m - r} \times \mathbf Z^{n - r}\) obtained from \(\lambda^\rho\) by deleting \(\lambda^\rho_{m_s}, \lambda^\rho_{\dot{n_s}}\) for \(1 \leq s \leq r\).

A weight is called \emph{typical} when \(r = 0\) and  \emph{\(r\)-fold atypical} or \emph{atypical} otherwise.
The matrix \(A(\lambda)\) with entries \(A(\lambda)_{ij} = (\lambda^\rho, \epsilon_i - \delta_j)\) is called the \emph{atypicality matrix} \cite{vhkt1990:character}.

There is a partial order on the set of weights such that  \(\mu \preceq \lambda\) if and only if \(\lambda\) and \(\mu\) have the same degree of atypicality, \(\typ{\mu} = \typ{\lambda}\), and \(\atyp{\mu} \leq \atyp{\lambda}\).
We use the convention that if \(a, b \in \mathbf Z^r\) then \(a \leq b\) if and only if \(a_i \leq b_i\) for all \(i\).
The \texttt{<=} operator has been overloaded to compare weights with respect to this partial order.

If \(\gamma_s\)is an atypical root of \(\lambda\) then the \emph{\(\gamma_s\)-height} of \(\lambda\) is the invariant
\begin{equation}
  h_s(\lambda) = \lambda_{m_s} - n_s + s.
\end{equation}
The \emph{height vector} and the \emph{height} of \(\lambda\) are the quantities
\begin{equation*}
  h(\lambda) = (h_1(\lambda), \dots, h_r(\lambda)),
  \qquad
  |h(\lambda)| = \sum_{s = 1}^r h_s(\lambda).
\end{equation*}

The invariants above are implemented as properties of the \texttt{Weight} object:
\begin{itemize}
\item \texttt{Weight.atypical\_roots}, \texttt{Weight.adeg}, and \texttt{Weight.atypicality\_matrix} calculate the list of atypical roots (as indices), the degree of atypicality, and the atypicality matrix respectively,
\item \texttt{Weight.typ} and \texttt{Weight.atyp} calculate the typical and the atypical tuple of a weight,
\item \texttt{Weight.height} computes the height vector.
\end{itemize}

The height of a weight with respect to its atypical roots is an important representation theoretic invariant.
The height vector completely determine the Kazhdan-Lusztig polynomials that determine the characters of irreducible \(\gl(m|n)\)-representations \cite{sz2007:character}.
The height vector is also used in constructing an equivalence of categories from a block \(\mathcal F^{m|n}_\chi\) of \(\gl(m|n)\) representations with central character \(\chi\) and atypicality \(r\) to the unique maximally atypical block \(\mathcal F^{r|r}_{\chi'}\) of \(\gl(r|r)\) \cite{ser1996:kazhdanlusztig, sz2007:character}.
Our package implements functionality to work with his equivalence of categories:
\begin{itemize}
  \item \texttt{height\_to\_atyp(height\_vector, typ)} determines the atypical tuple from a given height vector and typical tuple.
  \item \texttt{typ\_atyp\_to\_weight(typ, atyp)} calculates the weight in the block determined by the given typical tuple with the prescribed atypical tuple.
\end{itemize}

\begin{eg}
  \label{eg:rep-invariants}
  We follow the calculation in \cite[Example~3.3]{sz2007:character} to show the functionality of our package.
  As before, with
  \begin{equation*}
    \lambda = (7, 6, 5, 5, 3, 3, 2, 2, 0 \mid 1, 2, 3, 4, 4, 5, 7, 7) \in \mathbf Z^{9|8}
  \end{equation*}
  the atypical roots are given by
  \begin{align*}
    \Gamma(\lambda)
    &=  \left\{ \epsilon_8 - \epsilon_{11}, \epsilon_5 - \epsilon_{13}, \epsilon_{4} - \epsilon_{15}, \epsilon_2 - \epsilon_{16} \right\} \\
    &= \left\{ \epsilon_8 - \delta_2, \epsilon_5 - \delta_4, \epsilon_4 - \delta_6, \epsilon_2 - \delta_7 \right\}
  \end{align*}
  so \(\lambda\) is \(4\)-fold atypical.
  The equivalent calculation in Sage is as follows:
    \begin{lstlisting}
sage: from liesuperalg import Weight
sage: w = Weight([7, 6, 5, 5, 3, 3, 2, 2, 0], [1, 2, 3, 4, 4, 5, 7, 7])
sage: w.atypical_roots
[(8, 11), (5, 13), (4, 15), (2, 16)]
sage: w.adeg # degree of atypicality
4
sage: matrix(w.atypicality_matrix)
[ 14  12  10   8   7   5   2   1]
[ 12  10   8   6   5   3   0  -1]
[ 10   8   6   4   3   1  -2  -3]
[  9   7   5   3   2   0  -3  -4]
[  6   4   2   0  -1  -3  -6  -7]
[  5   3   1  -1  -2  -4  -7  -8]
[  3   1  -1  -3  -4  -6  -9 -10]
[  2   0  -2  -4  -5  -7 -10 -11]
[ -1  -3  -5  -7  -8 -10 -13 -14]
\end{lstlisting}

  The other invariants 
  \begin{align*}
    \atyp(\lambda) &=  (4, 8, 11, 14), \\
    \typ(\lambda) &= ((16, 12, 7, 5, 1), (2, 6, 9, 15)), \\
    h(\lambda) &= (1, 1, 2, 3)
  \end{align*}
  can be obtained by continuing the session as follows:
\begin{lstlisting}
sage: w.typ # typical tuple
[[16, 12, 7, 5, 1], [2, 6, 9, 15]]
sage: w.atyp # atypical tuple
[4, 8, 11, 14]
sage: w.height # height vector
[1, 1, 2, 3]
\end{lstlisting}

  Finally we check that the conversions related to equivalence of categories works as expected:
  \begin{lstlisting}
sage: from liesuperalg import typ_atyp_to_weight
sage: from liesuperalg import height_to_atyp
sage: t = [[16, 12, 7, 5, 1], [2, 6, 9, 15]]
sage: h = [1, 1, 2, 3]
sage: typ_atyp_to_weight(t, height_to_atyp(h, t))
gl(9|8) weight (7, 6, 5, 5, 3, 3, 2, 2, 0 | 1, 2, 3, 4, 4, 5, 7, 7)
sage: typ_atyp_to_weight(t, height_to_atyp(h, t)) == w
True
\end{lstlisting}
\end{eg}

\subsection{Cup diagrams}
We end this section with a description of cup diagrams associated to weights of \(\gl(m|n)\).
These diagrams are combinatorial gadgets that simplify many constructions in representation theory of lie superalgebras.
For a systematic introductions to these diagrams appears and their various applications in representation theory, see \cite{bs2010:highest-1, bs2010:highest-2, bs2010:highest-3, bs2010:highest-4}.

Let \(\lambda\) is an integral weight of \(\gl(m|n)\) and suppose \(\typ{\lambda} = (T^0, T^1)\).
The \emph{weight diagram} \(D(\lambda)\) of \(\lambda\) is  an assignment \(D(\lambda) \colon \mathbf Z \to \{ \varnothing , \bullet, \times, \blacktriangledown \}\) such that
\begin{equation*}
  D(\lambda)(k) =
  \begin{cases}
    \bullet & k \in T^0 \\
    \times & k \in T^1 \\
    \blacktriangledown & k \in \atyp{\lambda} \\
    \varnothing & \text{ otherwise}
  \end{cases}
\end{equation*}
Note that a weight diagram can be visualized on the number line \(\mathbf Z\).
The \emph{cup diagram} associated to the weight is obtained by starting, in order, from the left and joining a \(\blacktriangledown\) symbol with \(\varnothing\) such that all symbols appearing in between them are either \(\bullet\) or \(\times\) or are already joined by arcs.
We will often often use the symbol \(\blacktriangle\) for the \(\varnothing\) matching with \(\blacktriangledown\).

For any given weight, the \texttt{cup\_diagram(weight)} produces this diagram.

\begin{eg}
  If \(\lambda\) is the weight from Example~\ref{eg:rep-invariants}, continuing the session with 
\begin{lstlisting}
sage: from liesuperalg import cup_diagram
sage: cup_diagram(w)
\end{lstlisting}
  produces the cup diagram
  \begin{equation*}
    \includegraphics[scale=0.7, trim=0 0 0.5in 1in, clip]{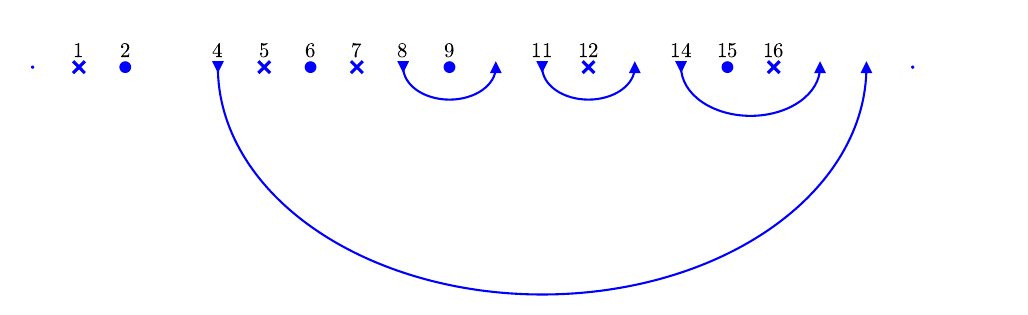}
  \end{equation*}
\end{eg}

\section{Representations, Characters, and Kazhdan-Lusztig Polynomials}
\label{sec:reps}
If \(V'\) is an ordinary (not graded) complex vector space, then the representation theory of \(\gl(V')\) is described by highest weight theory \cite{FH:RepresentationTheory, Hum:IntroductionLie}.
We denote the representation of \(\gl(V')\) corresponding to the weight \(\mu\) by \(\SS_\mu(V')\).
Highest weight theory has been generalized to the case of \(\gl(m|n)\) and we describe the basic notions here.
Since \(\gl(m|n)\) contains as subalgebra \(\g_0 = \gl(m) \oplus \gl(n)\), a class of representations of \(\gl(m|n)\) can be constructed by inducing from representations
\begin{equation*}
  L_0(\lambda) = \SS_{\lambda^0} V_0 \otimes \SS_{\lambda^1} V_1
\end{equation*}
of \(\g_0\).
If \(\lambda\) is a dominant integral weight of \(\gl(m|n)\), we define the \emph{Kac module} as the induced module
\begin{equation}
  \label{eq:def:kac-module}
  K(\lambda) =
  \mathrm{Ind}^\g_{\gl(m) \oplus \gl(n)} L_0(\lambda) =
  \mathscr U(V_0^\ast \otimes V_1) \otimes L_0(\lambda).
\end{equation}
where \(\lambda^0 = (\lambda_1, \dots, \lambda_m)\) and \(\lambda^1 = (\lambda_{\dot{1}}, \dots, \lambda_{\dot{n}})\) are weights of \(\gl(m)\) and \(\gl(n)\) respectively and \(\mathscr U(V_0^\ast \otimes V_1)\) is the universal enveloping algebra of \(V_0^\ast \otimes V_1\).
Kac modules are indecomposable but, in general, not irreducible.

\begin{thm}[{\cite{kac1977:lie}}]
  The module \(K(\lambda)\) has a unique quotient \(L(\lambda)\). If \(\lambda\) is dominant and integral then \(L(\lambda)\) is irreducible. Furthermore, all irreducible \(\g\)-modules are obtained in this way.
\end{thm}

The problem of determining a character formula for \(L(\lambda)\) and the composition factors of \(K(\lambda)\) has been a motivating problem in representation theory of \(\gl(m|n)\).
In particular, one is generally interested in determining the character \(\ch{L(\lambda)}\) or \(\ch{K(\lambda)}\) in terms of a suitable character basis by determining, for example, the coefficients \(a_{\lambda\mu}\) or \(b_{\lambda\mu}\) such that 
\begin{equation*}
  \ch{K(\lambda)} = \sum_\mu a_{\lambda\mu} \ch{L(\mu)},
  \qquad
  \ch{L(\lambda)} = \sum_\mu b_{\lambda\mu} \ch{K(\mu)}.
\end{equation*}
In \cite{ser1998:characters}, Serganova determined a character formula for \(L(\lambda)\) in terms of generalized Kazhdan-Lusztig polynomials \(K_{\lambda,\mu}(q)\) and showed that the coefficients \(b_{\lambda\mu}\) above are given by
\begin{equation*}
  b_{\lambda\mu} = K_{\lambda,\mu}(-1).
\end{equation*}
In \cite{bru2003:kazhdanlusztig}, Brundan provided an algorithm for calculating these polynomials and determining  the composition factors of \(K(\lambda)\).
Our main reference will be the explicit calculation of the character formula proven in \cite{sz2007:character} using Brundan's algorithm.

If \(\lambda\) and \(\mu\) satisfy \(\mu \preceq \lambda\) then they have the same degree of atypicality; say \(\lambda\) and \(\mu\) are both \(r\)-fold atypical.
Following the form of Brundan's algorithm as described in \cite{sz2007:character}, we first calculate the \emph{distance} between two atypical roots \(\gamma_s, \gamma_t\) as
\begin{equation*}
  d_{s, t}(\lambda) = \Bigl| [a_s, a_t] \setminus D(\lambda) ^{-1} \left( \{ \bullet, \times, \blacktriangledown \} \right) \Bigr|
\end{equation*}
where \(a_s\) and \(a_t\) are the \(s\)-th and \(t\)-th entries of \(\atyp{\lambda}\).
Two atypical roots are \(c\)-related if \(d_{s, t}(\lambda) < t - s\) and are strongly \(c\)-related if \(\gamma_i\) and \(\gamma_{i + 1}\) are \(c\)-related for all \(s \leq i < t\).
Since the permutations \(\mathrm{Sym}(\Gamma(\lambda)) \simeq \mathfrak S_r\) form a subgroup of the Weyl group of \(\gl(m|n)\), the rough idea is to understand the dot action of these elements on \(\lambda\) satisfying \(\mu \preceq \sigma \bullet \lambda\) such that \(\sigma\) does not change the order of \(s < t\) whenever \(\gamma_s\) and \(\gamma_t\) are strongly \(c\)-related.
If \(S^{\lambda, \mu}\) is the set of all such permutations, then following theorem calculates the generalized Kazhdan-Lusztig polynomial \(K_{\lambda, \mu}(q)\):
\begin{thm}[{\cite[Theorem~3.23]{sz2007:character}}]
  Suppose \(\lambda\) and \(\mu\) are dominant weights. Then \(K_{\lambda, \mu}(q) \neq 0\) if and only if \(\mu \preceq \lambda\) and in that case
  \begin{equation*}
    K_{\lambda, \mu}(q) = q^{|h(\lambda)| - |h(\mu)|} \sum_{\sigma \in S^{\lambda, \mu}} q^{-2l(\sigma)}
  \end{equation*}
  where \(l(\sigma)\) is Bruhat length of \(\sigma\).
\end{thm}

The details of this construction are quite technical so we refer reader to \cite{sz2007:character}.
If necessary the \texttt{Weight} object exposes some of this technical machinery:
\begin{itemize}
\item \texttt{Weight.cr} and \texttt{Weight.scr} produce dictionaries encoding \(c\)-connectivity and strong \(c\)-connectivity of atypical roots.
\item \texttt{Weight.atyp\_dot\_action(self, sigma)} applies the dot action of elements of \(\sigma \in \mathrm{Sym}(\Gamma(\lambda))\) to \(\lambda\).
\item \texttt{Weight.respects\_scr(self, sigma)} determines whether a given permutation \(\sigma\) does not change the order of \(s < t\) for strongly \(c\)-connected \(\gamma_s\) and \(\gamma_t\).
\end{itemize}

The package exposes the two main methods to work with characters of \(\gl(m|n)\):
\begin{itemize}
\item \texttt{gen\_KL(weight\_lambda, weight\_mu, q)} calculates the generalized Kazhdan-Lusztig polynomial \(K_{\lambda, \mu}(q)\).
\item \texttt{mult\_kac\_in\_irrd(weight\_lambda, weight\_mu)} calculates the multiplicity of a Kac module in the expansion of the given irreducible.
\end{itemize}
\begin{eg}
  Suppose
  \begin{align*}
    \lambda
    &=  ( 7,  6,  5,  5,  3,  3,  2,  2,  0 \mid  1,  2,  3,  4,  4,  5,  7,  7)
    \\
    \mu
    &= ( 7,  4,  4,  4,  2,  1,  1,  1,  0 \mid  1,  1,  1,  2,  4,  4,  4,  7)
  \end{align*}
  then using Theorem~\cite[Theorem~3.23]{sz2007:character}, we obtain
  \begin{equation*}
    K_{\lambda, \mu}(q) = q^5\left(\frac{1}{q^2} + 1\right)
  \end{equation*}
  and so the multiplicity of \(\ch{K(\mu)}\) is \(K_{\lambda, \mu}(-1) = -2\).
  \begin{lstlisting}
sage: from liesuperalg import gen_KL
sage: from liesuperalg import mult_kac_in_irrd
sage: var('q')
sage: w_lambda = Weight([7, 6, 5, 5, 3, 3, 2, 2, 0], [1, 2, 3, 4, 4, 5, 7, 7])
sage: w_mu = Weight([7, 4, 4, 4, 2, 1, 1, 1, 0], [1, 1, 1, 2, 4, 4, 4, 7])
sage: w_lambda >= w_mu
True
sage: gen_KL(w_lambda, w_mu, q)
q^5*(1/q^2 + 1)
sage: mult_kac_in_irrd(w_lambda, w_mu)
-2
sage: mult_kac_in_irrd(w_lambda, w_mu) == gen_KL(w_lambda, w_mu, -1)
True
\end{lstlisting}
\end{eg}

\section{Composition Factors of Kac Modules}
\label{sec:composition}

Recall that if \(\lambda\) is a dominant weight then the Kac module \(K(\lambda)\) is indecomposable but not irreducible in general.
As such, one is interested in calculating the irreducible composition factors of the Kac module \(K(\lambda)\).
Brundan proves in \cite[Main Theorem]{bru2003:kazhdanlusztig} that all composition factors can be obtained by iterative applications of the so-called \emph{lowering operator} on the weight diagram associated to \(\lambda\). 
Our package implements this calculation based on the version  of the algorithm described in \cite{sz2012:generalised} and \cite[Theorem~1.1]{su2006:composition}
We refer the reader to these two sources for the numerous technical details involved in the construction of these operators.

The function \texttt{kac\_composition\_factors(weight\_lambda)} produces an iterator for all the composition factors \(K(\lambda)\).
Its usage is describe in the following example.

\begin{eg}
  If \(\lambda = (8, 5, 5, 3, 3, 2, 2 \mid 2, 3, 4, 4, 5, 9)\) then, as computed in \cite[Example~3.13]{su2006:composition}, the \(\rho\)-translates of the first \(7\) of the \(14\) composition factors are given by:
 \begin{align*}
   (15, 11, 10, 7, 6, 4, 3 &\mid 3, 5, 7, 8, 10, 15), \\
   (15, 11, 10, 7, 6, 4, 2 &\mid 2, 5, 7, 8, 10, 15), \\
   (15, 11, 10, 6, 4, 2, 1 &\mid 1, 2, 5, 8, 10, 15), \\
   (15, 11, 9, 7, 6, 4, 3 &\mid 3, 5, 7, 8, 9, 15),\\
   (15, 11, 9, 7, 6, 4, 2 &\mid 2, 5, 7, 8, 9, 15),\\
   (15, 11, 9, 6, 4, 2, 1 &\mid 1, 2, 5, 8, 9, 15),\\
   (15, 11, 7, 6, 4, 2, 1 &\mid 1, 2, 5, 7, 8, 15).
 \end{align*}

 We can compute the full set as follows:
  \begin{lstlisting}
sage: # Example 3.13 (Su 2006)
sage: from liesuperalg import kac_composition_factors
sage: w_lambda = Weight([8, 5, 5, 3, 3, 2, 2], [2, 3, 4, 4, 5, 9])
sage: for i, w_mu in enumerate(kac_composition_factors(w_lambda), start=1):
....:    print(f"{i:2}: {w_mu.rho()}")
....:
 1: (15, 11, 10,  7,  6,  4,  3 |  3,  5,  7,  8, 10, 15)
 2: (15, 11, 10,  7,  6,  4,  2 |  2,  5,  7,  8, 10, 15)
 3: (15, 11, 10,  6,  4,  2,  1 |  1,  2,  5,  8, 10, 15)
 4: (15, 11,  9,  7,  6,  4,  3 |  3,  5,  7,  8,  9, 15)
 5: (15, 11,  9,  7,  6,  4,  2 |  2,  5,  7,  8,  9, 15)
 6: (15, 11,  9,  6,  4,  2,  1 |  1,  2,  5,  8,  9, 15)
 7: (15, 11,  7,  6,  4,  2,  1 |  1,  2,  5,  7,  8, 15)
 8: (14, 11, 10,  7,  6,  4,  3 |  3,  5,  7,  8, 10, 14)
 9: (14, 11, 10,  7,  6,  4,  2 |  2,  5,  7,  8, 10, 14)
10: (14, 11, 10,  6,  4,  2,  1 |  1,  2,  5,  8, 10, 14)
11: (14, 11,  9,  7,  6,  4,  3 |  3,  5,  7,  8,  9, 14)
12: (14, 11,  9,  7,  6,  4,  2 |  2,  5,  7,  8,  9, 14)
13: (14, 11,  9,  6,  4,  2,  1 |  1,  2,  5,  8,  9, 14)
14: (14, 11,  7,  6,  4,  2,  1 |  1,  2,  5,  7,  8, 14)
\end{lstlisting}
\end{eg}
  
\begin{eg}
  If \(\lambda = (2, 1, 1, 0, 0 \mid 0, 0, 1, 3, 3, 4)\) then, as determined in \cite[Example~5.6]{sz2012:generalised}, the Kac module \(K(\lambda)\) has \(19\) composition factors.
  We can explicitly calculate these as follows:
  \begin{lstlisting}
sage: # Example 5.6 (Su-Zhang 2012)
sage: from liesuperalg import kac_composition_factors
sage: w_lambda = Weight([2, 1, 1, 0, 0], [0, 0, 1, 3, 3, 4])
sage: for i, w_mu in enumerate(kac_composition_factors(w_lambda), start=1):
....:     print(f"{i:2}: {w_mu}")    
....:
 1: ( 2,  1,  1,  0,  0 |  0,  0,  1,  3,  3,  4)
 2: ( 2,  1,  1,  0, -1 | -1,  0,  1,  3,  3,  4)
 3: ( 2,  1,  1, -2, -2 | -2, -2,  1,  3,  3,  4)
 4: ( 2,  1,  0,  0,  0 |  0,  0,  0,  3,  3,  4)
 5: ( 2,  1,  0,  0, -1 | -1,  0,  0,  3,  3,  4)
 6: ( 2,  1,  0, -2, -2 | -2, -2,  0,  3,  3,  4)
 7: ( 2,  1, -1, -2, -2 | -2, -2, -1,  3,  3,  4)
 8: ( 1,  1,  1,  0,  0 |  0,  0,  1,  2,  3,  4)
 9: ( 1,  1,  1,  0, -1 | -1,  0,  1,  2,  3,  4)
10: ( 1,  1,  1, -2, -2 | -2, -2,  1,  2,  3,  4)
11: ( 1,  1,  0,  0,  0 |  0,  0,  0,  2,  3,  4)
12: ( 1,  1,  0,  0, -1 | -1,  0,  0,  2,  3,  4)
13: ( 1,  1,  0, -2, -2 | -2, -2,  0,  2,  3,  4)
14: ( 1,  1, -1, -2, -2 | -2, -2, -1,  2,  3,  4)
15: ( 0,  0,  0,  0,  0 |  0,  0,  0,  0,  3,  4)
16: ( 0,  0,  0,  0, -1 | -1,  0,  0,  0,  3,  4)
17: ( 0,  0,  0, -2, -2 | -2, -2,  0,  0,  3,  4)
18: ( 0,  0, -1, -2, -2 | -2, -2, -1,  0,  3,  4)
19: ( 0, -1, -1, -2, -2 | -2, -2, -1, -1,  3,  4)
\end{lstlisting}
\end{eg}

\section{Decomposition of Lie Superalgebra Modules}
\label{sec:decomposition}
Finally we describe an algorithm to decompose characters of \(\gl(m|n)\)-modules in terms of the irreducible characters.
Our main goal is to determine whether a given a \(\g_0\)-module
\begin{equation*}
  M = \bigoplus_\alpha L_0(\alpha)^{\oplus m_\alpha}
\end{equation*}
can be written as a sum of irreducible modules in the Grothendieck group of \(\gl(m|n)\)-modules.

Since
\begin{equation*}
  K(\mu) = \mathscr U(V_0^\ast \otimes V_1) \otimes L_0(\mu) = \bigwedge (V_0^\ast \otimes V_1) \otimes L_0(\mu), 
\end{equation*}
we may use classical representation theory, in particular the Littlewood-Richardson rule, to determine the character decomposition
\begin{equation*}
  \ch K(\mu) = \sum_\alpha c_{\mu, \alpha} \ch{L_0(\alpha)}
\end{equation*}
where \(c_{\mu, \alpha}\) is a sum of products of Littlewood-Richardson coefficients.
The characters of the irreducibles are given by
\begin{equation*}
  \ch{L(\lambda)}
  = \sum_\mu K_{\lambda, \mu}(-1) \ch{K(\mu)}
  = \sum_{\mu, \alpha} K_{\lambda, \mu}(-1) \cdot c_{\mu, \alpha} \ch{L_0(\alpha)}.
\end{equation*}

Thus, given \(\ch{M} = \sum_\alpha m_\alpha L_0(\alpha)\) our algorithm determines a set of \(\lambda\) and multiplicities \(n_\lambda\) such that 
\begin{equation*}
  m_\alpha = \sum n_\lambda \cdot K_{\lambda, \mu}(-1) \cdot c_{\mu, \alpha}.
\end{equation*}

This is done recursively in the obvious way:
\begin{enumerate}
\item If all \(m_\alpha = 0\), a decomposition is given by \(0\). In this case, terminate.
\item If there exist non-zero \(m_\alpha\):
  \begin{enumerate}
  \item 
    Find the largest weight \(\alpha\) with \(m_\alpha \neq 0\) and determine its \emph{support} consisting of Kac modules.
    That is find all \(\mu\) such that \(c_{\mu, \alpha} \neq 0\).
  \item Since \(K_{\alpha, \alpha}(-1) \cdot c_{\alpha, \alpha} = 1\), pick \(\ch{L(\alpha)}\) as a summand in the decomposition.
  \item Recursively decompose
    \begin{equation*}
      \ch{M} - \sum_\beta \sum_{\mu \in \mathrm{supp(\alpha)}} K_{\alpha, \mu}(-1) \cdot c_{\mu, \beta} \ch{L_0(\beta)}
    \end{equation*}
  \item Add \(\ch{L(\alpha)}\) to the decomposition of the smaller module and terminate.
  \end{enumerate}
\end{enumerate}

The command that implements this procedure is \texttt{decompose(module, cache=None)}  accepts a module given as a dictionary \(\{ \alpha \colon m_\alpha \}\).
Since computing the Kac module support of any given \(L_0(\alpha)\) is computationally very expensive, a pre-computed cache can be specified by setting the optional \texttt{cache} parameter.
The return value is a dictionary \(\{ \lambda \colon n_\lambda \}\) specifying the decomposition.

\begin{eg}
  The following sage session calculates the decomposition of the module \(M\) given below.
\begin{lstlisting}
sage: from liesuperalg import decompose
sage: M = {
....:     Weight([ 1, -2, -2, -2, -2, -2], [-2, -2, -2, -2, -2, -1]) : 1,  
....:     Weight([ 1, -2, -2, -2, -2, -3], [-2, -2, -2, -2, -2, -2]) : 1,   
....:     Weight([ 0, -2, -2, -2, -2, -2], [-2, -2, -2, -2, -2, -2]) : 1,   
....:     Weight([ 0, -2, -2, -2, -2, -2], [-3, -2, -2, -2, -2, -1]) : 1,   
....:     Weight([ 0, -2, -2, -2, -2, -3], [-3, -2, -2, -2, -2, -2]) : 1,   
....:     Weight([-1, -2, -2, -2, -2, -2], [-3, -2, -2, -2, -2, -2]) : 1,   
....:     Weight([-1, -2, -2, -2, -2, -2], [-3, -3, -2, -2, -2, -1]) : 1,   
....:     Weight([-1, -2, -2, -2, -2, -3], [-3, -3, -2, -2, -2, -2]) : 1,   
....:     Weight([-2, -2, -2, -2, -2, -2], [-3, -3, -2, -2, -2, -2]) : 1,   
....:     Weight([-2, -2, -2, -2, -2, -2], [-3, -3, -3, -2, -2, -1]) : 1,   
....:     Weight([-2, -2, -2, -2, -2, -3], [-3, -3, -3, -2, -2, -2]) : 1
....: }
sage: decompose(M)
{gl(6 | 6) weight ( 1, -2, -2, -2, -2, -2 | -2, -2, -2, -2, -2, -1): 1}
\end{lstlisting}  
\end{eg}

\bibliographystyle{alpha}
\bibliography{refs}
\end{document}